\documentclass[12pt,a4paper]{amsart}

\usepackage{color}
\usepackage{amscd,amssymb,amsopn,amsmath,amsthm,
                   graphics,amsfonts,enumerate,verbatim,calc,mathdots,mathrsfs}
\usepackage[dvips]{graphicx}
\usepackage{tikz}
\usetikzlibrary{intersections, calc, arrows}

\newtheorem{Theorem}{Theorem}[section]
\newtheorem{Proposition}[Theorem]{Proposition}
\newtheorem{Lemma}[Theorem]{Lemma}

\newtheorem{Remark}[Theorem]{Remark}
\newtheorem{Example}[Theorem]{Example}
\newtheorem{Definition}[Theorem]{Definition}

\newtheorem{Claim}[Theorem]{Claim}

\makeatletter
  
  \@addtoreset{equation}{section}
\makeatother

\newcommand{\gp}{\mathfrak{p}}

\newcommand{\gm}{\mathfrak{m}}

\newcommand{\ga}{\mathfrak{a}}

\newcommand{\pp}{\mathbb{P}}
\newcommand{\zz}{\mathbb{Z}}
\newcommand{\nn}{\mathbb{N}}

\newcommand{\calP}{\mathscr{P}}

\newcommand{\rees}[1]{{\mathscr R}(#1)}
\newcommand{\gr}[1]{{\mathscr G}(#1)}

\newcommand{\srees}[1]{{\mathscr R}_s(#1\,)}
\newcommand{\esrees}[1]{{\mathscr R}'_s(#1)}

\newcommand{\mult}[2]{{\mathrm e}_{#1}(#2)} %multiplicity
\newcommand{\length}[2]{\ell_{#1}(\,#2\,)}

\newcommand{\Min}[2]{{\rm Min}_{#1}\,{#2}}
\newcommand{\height}[2]{{\rm ht}_{#1}\,#2}

\begin{document}
\title[Finitely generated symbolic Rees rings]{
Finitely generated symbolic Rees rings of ideals defining
certain finite sets of points in $\pp^2$}

\author[K. Kai]{Keisuke Kai}
\address[K. Kai]{
Department of Mathematics and Informatics,
Graduate School of Science,
Chiba University,
Yayoi-cho 1-33, Inage-ku,
Chiba 263-8522, Japan}
\email{k\_k04\_1023@yahoo.co.jp}

\author[K. Nishida]{Koji Nishida}
\address[K. Nishida, corresponding author]{
Department of Mathematics and Informatics,
Graduate School of Science,
Chiba University,
Yayoi-cho 1-33, Inage-ku,
Chiba 263-8522, Japan}
\email{nishida@math.s.chiba-u.ac.jp}

\maketitle

\begin{abstract}
The purpose of this paper is to prove that
the symbolic Rees rings of ideals defining certain finite sets
of points in the projective plane over an algebraically closed field
are finitely generated using a ring theoretical criterion
which is known as Huneke's criterion.
\end{abstract}

\section{Introduction}

Let $R$ be a commutative Noetherian ring and
let $\ga$ be a proper ideal of $R$.
We denote the set of minimal prime divisors of $\ga$ by $\Min{}{\ga}$.
For any $r \in \zz$, we define
\[
\ga^{(r)} = \bigcap_{\gp\,\in\,\Min{}{\ga}} (\gp^rR_\gp \cap R)
\]
and call it the $r$-th symbolic power of $\ga$.
Moreover, taking an indeterminate $t$, 
we define the symbolic Rees ring of $\ga$ by
\[
\srees{\ga} = \sum_{r \in \nn} \ga^{(r)}t^r \subset R[ t ]\,,
\]
where $\nn = \{ 0, 1, 2, \dots \}$.
Although deciding whether the symbolic Rees rings of given ideals
are finitely generated or not is an important problem
in commutative algebra and algebraic geometry,
but usually it is a hard task.
In this paper, we focus our attention on a ring theoretical criterion
for finite generation of symbolic Rees rings
which is known as Huneke's criterion in a special situation described below.

Let $K$ be a field and
let $I$ be a proper homogeneous ideal of the polynomial ring
$S = K[x, y, z]$ which we regard as an $\nn$-graded ring by setting
the degrees of $x$, $y$ and $z$ to suitable positive integers.
We assume that $S / I$ is a $1$-dimensional reduced ring.
Let $\gm = (x, y, z)S$.
Because the symbolic powers of $I$ are also homogeneous,
we have $S_\gm \otimes_S \srees{I} = \srees{IS_\gm}$,
i.e., $I^{(r)}S_\gm = (IS_\gm)^{(r)}$ for any $r \in \zz$.
On the other hand,
if $\gp \in \Min{}{I}$,
we have $IS_\gp = \gp S_\gp$ as $\sqrt{I} = I$,
and so $S_\gp \otimes_S \srees{I}$ coincides with
\[
\rees{S_\gp} = \sum_{r \in \nn} \,\gp^rS_\gp\cdot t^r\,,
\]
which is the ordinary Rees ring of the $2$-dimensional regular local ring $S_\gp$.
Here, let us recall the following condition introduced in
\cite[Theorem 3.25]{H} and \cite[Proposition 2.1]{KN}.
\begin{Definition}\label{1.1}
Let $0 < r_i \in \nn$ and $\xi_i \in I^{(r_i)}$ for $i = 1, 2$.
We say that $\xi_1$ and $\xi_2$ satisfy Huneke's condition on $I$
if the following two equalities hold.
\begin{itemize}
\item[{\rm (a)}]
$IS_\gm = \sqrt{(\xi_1, \xi_2)S_\gm}$\,.
\item[{\rm (b)}]
$\gr{S_\gp}_+ = \sqrt{(\xi_1t^{r_1}, \xi_2t^{r_2})\gr{S_\gp}}$ for any $\gp \in \Min{}{I}$\,,
\end{itemize}
where $\gr{S_\gp} = S_\gp / \gp S_\gp \otimes \rees{S_\gp}$
and $\gr{S_\gp}_+$ denotes the ideal generated by the homogeneous elements
of positive degree.
In both of these equalities,
the right sides are obviously contained in the left sides,
so the crucial requirement of the condition stated above is that
the left sides are included in the right sides.
\end{Definition}
Although the condition stated in Definition \ref{1.1} is rather complicated,
it is equivalent to an easy condition if 
the grading of $S$ is ordinary and
both of $\xi_1$ and $\xi_2$ are homogeneous.
\begin{Proposition}\label{1.2}
Suppose $\deg x = \deg y = \deg z = 1$.
Let $0 < r_i, d_i \in \nn$ and
$\xi_i \in [I^{(r_i)}]_{d_i}$ for $i = 1, 2$.
Then $\xi_1$ and $\xi_2$ satisfy Huneke's condition on $I$
if and only if $\height{}{(\xi_1, \xi_2)S} = 2$ and
\[
\frac{d_1}{r_1}\cdot\frac{d_2}{r_2} = \mult{}{S / I}\,,
\]
where $\mult{}{S / I}$ denotes the multiplicity of $S / I$
{\rm (cf. \cite[Definition 4.1.5]{BH})} .
\end{Proposition}
Now, Huneke's criterion can be described as follows.
\begin{Theorem}\label{1.3}
$\srees{I}$ is finitely generated if and only if
there exist elements in $I^{(r_1)}$ and $I^{(r_2)}$ satisfying
Huneke's condition on $I$ for some $0 < r_1, r_2 \in \nn$.
\end{Theorem}
Huneke's criterion was originally proved by Huneke
(cf. \cite[Theorem 3.1, 3.2]{H})
in the case where $I$ is a prime ideal,
and the generalized version was given by Kurano and Nishida
(cf. \cite[Theorem 2.5]{KN}) 
so that it can be applied to radical ideals.
The purpose of this paper is to prove that
the symbolic Rees rings of the ideals defining certain finite sets
in the projective plane $\pp^2$ are finitely generated using Huneke's criterion.

Let $K$ be an algebraically closed field
and $\deg x = \deg y = \deg z = 1$.
For a point $P = (a : b : c) \in \pp^2 = \pp_K^{\,2}$,
we denote by $I_P$ the ideal of $S$ generated by the maximal minors of the matrix
\[
\left(\begin{array}{ccc}
x & y & z \\
a & b & c
\end{array}\right)\,,
\]
which is the defining ideal of $P$.
Of course, $I_P$ is a prime ideal of $S$ generated by a regular sequence.
Moreover, for a set $H = \{ P_1, P_2, \dots , P_e \}$
of $e$ points in $\pp^2$, we set
\[
I_H = I_{P_1} \cap I_{P_2} \cap \cdots \cap I_{P_e}\,.
\]
Then we have
\[
I_H^{\,(r)} = I_{P_1}^{\,r} \cap I_{P_2}^{\,r} \cap \cdots \cap I_{P_e}^{\,r}
\]
for any $r \in \zz$.
As is well known, $\srees{I_H}$ is finitely generated if and only if so is
\[
\esrees{I_H} = \sum_{r \in \zz}\,I_H^{\,(r)}t^r \subset S[ t, t^{-1} ]\,,
\]
and the finite generation of these graded rings is related to
that of the Cox ring $\Delta_H$, which is the subring
\[
\sum_{(r_1, \dots , r_e) \in \zz^e}
(I_{P_1}^{\,r_1} \cap \cdots \cap I_{P_e}^{\,r_e})\,t_1^{\,r_1} \cdots t_e^{\,r_e}
\]
of $S[t_1^{\,\pm 1}, \dots , t_e^{\,\pm 1}]$,
where $t_1, \dots , t_e$ are indeterminates.
Since $\esrees{I_H}$ coincides with the diagonal part of $\Delta_H$,
$\esrees{I_H}$ is finitely generated if so is $\Delta_H$.
For example, in \cite{EKW}
Elizondo, Kurano and Watanabe proved that
$\Delta_H$ is finitely generated if the points of $H$ lie on a line in $\pp^2$.
Moreover, in \cite{TVV} Testa, Varilly-Alvarado and Velasco proved
the finite generation of $\Delta_H$ for the following cases.
\begin{itemize}
\item[(i)]
$e \leq 8$.
\item[(ii)]
$e - 1$ points in $H$ lie on a (possibly reducible) conic in $\pp^2$.
\item[(iii)]
$H$ consists of $10$ points of pairwise intersections of
$5$ general lines in $\pp^2$.
\item[(iv)]
There exist $3$ distinct lines $L_1$, $L_2$ and $L_3$ in $\pp^2$
such that $H$ consists of pairwise intersections of these lines
and $2$, $3$ and $5$ additional points on
$L_1$, $L_2$ and $L_3$, respectively ($e = 13$).
\end{itemize} 
Of course, $\srees{I_H}$ can be finitely generated for wider classes of $H$.
For example, the following is known. 
\begin{Theorem}\label{1.4}
Let $n$ be a positive integer which is not a multiple
of the characteristic of $K$ and let $\theta$ be a primitive $n$-th root of unity.
We set
\[
H = \{ (1 : 0 : 0) , (0 : 1 : 0) , (0 : 0 : 1) \} \cup
\{ (\theta^i : \theta^j : 1) \mid i , j = 1, \dots , n \}\,.
\]
Then $\srees{I_H}$ is finitely generated.
\end{Theorem}
If $n = 1$ or $2$, then the number of points in $H$ stated
in the above theorem is $4$ or $7$,
and so the finite generation of $\srees{I_H}$ follows from that of $\Delta_H$.
In \cite{HS}, Harbourne and Seceleanu proved Theorem \ref{1.4}
in the case where $n = 3$,
and the case where $n \geq 4$ was settled by Nagel and Seceleanu in \cite{NS}.
In this paper, we aim to give an alternative proof for Theorem \ref{1.4}
using Huneke's criterion.
In Section 3,
we will show that there exist two elements in $I_H^{\,(n)}$
satisfying Huneke's condition on $I_H$.
Although both of those elements are homogeneous in the case where $n = 3$,
but one of the two elements is not homogeneous if $n \geq 4$.
Moreover, by a similar argument we prove that
the following assertion holds.
\begin{Theorem}\label{1.5}
Let $f$ and $g$ be homogeneous polynomials in $S$
such that $S / (f, g)$ is a $1$-dimensional reduced ring.
We put $\deg f = m$ and $\deg g = n$.
Let us assume that
\[
f \in I_A^{\,m}\,,\hspace{1ex}
g \in I_B^{\,n}\,,\hspace{1ex}
f \not\in I_B\hspace{1ex} \mbox{and} \hspace{1ex}
g \not\in I_A\,,
\]
where $A$ and $B$ are distinct two points  in $\pp^2$.
We set
\[
H = \{ A\,,\, B \} \cup
\{ P \in \pp^2 \mid (f , g) \subseteq I_P \}\,.
\]
Then $\srees{I_H}$ is finitely generated.
\end{Theorem}
The above theorem will be proved in Section 4 showing that
there exist linear forms $f_1, f_2, \dots, f_m \in [ I_A ]_1$ and
$g_1, g_2, \dots, g_n \in [ I_B ]_1$ such that
\begin{eqnarray*}
 & & f = f_1f_2 \cdots f_m \,\,, \hspace{1ex}
 g = g_1g_2 \cdots g_n\,\,, \\
 & & \mbox{$f_i \not\in I_B$ for any $i = 1, 2, \dots , m$ \,and} \\
 & & \mbox{$g_j \not\in I_A$ for any $j = 1, 2, \dots , n$\,.} 
\end{eqnarray*}
Let $P_{ij}$ be the intersection point of the lines defined by $f_i$ and $g_j$.
Because $S / (f, g)$ is reduced,
$f_i \not\sim f_k$
(i.e., $f_i / f_k \not\in K$)  if $i \neq k$,
and $g_j \not\sim g_\ell$ if $j \neq \ell$.
Consequently, we see
\[
H = \{ A, B \} \cup \{ P_{ij} \mid 
\mbox{$i = 1, \dots, m$ and $j = 1, \dots, n$} \}
\]
and $\sharp H = mn + 2$ (Figure \ref{Fig:1}).
We will prove that $\srees{I_H}$ is finitely generated
by finding elements in
$I_H^{\,(mn)}$ and $I_H^{\,(2)}$ satisfying Huneke's condition on $I_H$.
If $m \neq n$, then both of those elements are not homogeneous.

\begin{figure}[h]
\begin{tikzpicture}
\draw[thick, name path=f1] (0,2.75)--(1,3)--(5,4);
\draw[thick, name path=f2] (0,3)--(1,3)--(5,3);
\draw[thick, dotted] (5.3,2.5)--(5.3,2);
\draw[thick, name path=fm] (0,3.375)--(1,3)--(5,1.5);
\draw[thick, name path=g1] (3.25,6)--(3,5)--(2,1);
\draw[thick, name path=g2] (3,6)--(3,5)--(3,1);
\draw[thick, dotted] (3.5,0.75)--(4,0.75);
\draw[thick, name path=gn] (2.625,6)--(3,5)--(4.5,1);
\coordinate (A) at (1,3) node at (A) [below] {$A$};
\coordinate (B) at (3,5) node at (B) [right] {$B$};
\node (f1) at (5,4) [right] {$f_1$};
\node (f2) at (5,3) [right] {$f_2$};
\node (fm) at (5,1.5) [right] {$f_m$};
\node (g1) at (2,1) [below] {$g_1$};
\node (g2) at (3,1) [below] {$g_2$};
\node (gn) at (4.5,1) [below] {$g_n$};
\fill (A) circle (2pt) (B) circle (2pt);
\path [name intersections={of=f1 and g1}];
\coordinate (P11) at (intersection-1);
\fill (P11) circle (2pt);
\path [name intersections={of=f1 and g2}];
\coordinate (P12) at (intersection-1);
\fill (P12) circle (2pt);
\path [name intersections={of=f1 and gn}];
\coordinate (P1n) at (intersection-1);
\fill (P1n) circle (2pt);
\path [name intersections={of=f2 and g1}];
\coordinate (P21) at (intersection-1);
\fill (P21) circle (2pt);
\path [name intersections={of=f2 and g2}];
\coordinate (P22) at (intersection-1);
\fill (P22) circle (2pt);
\path [name intersections={of=f2 and gn}];
\coordinate (P2n) at (intersection-1);
\fill (P2n) circle (2pt);
\path [name intersections={of=fm and g1}];
\coordinate (Pm1) at (intersection-1);
\fill (Pm1) circle (2pt);
\path [name intersections={of=fm and g2}];
\coordinate (Pm2) at (intersection-1);
\fill (Pm2) circle (2pt);
\path [name intersections={of=fm and gn}];
\coordinate (Pmn) at (intersection-1);
\fill (Pmn) circle (2pt);
\end{tikzpicture}
\caption{Theorem \ref{1.5}}
\label{Fig:1}
\end{figure}

Setting $f = y^m - z^m$, $g = z^n - x^n$, $A = (1 : 0 : 0)$
and $(0 : 1 : 0)$ in Theorem \ref{1.5},
we get the following example.

\begin{Example}\label{1.6}
Let $m, n$ be positive integers which are not multiples of the
characteristic of $K$.
Let $\theta_m$ and $\theta_n$ be primitive 
$m$-th and $n$-th root of unity, respectively.
We set
\[
H = \{ (1 : 0 : 0), (0 : 1 : 0) \} \cup
\{ (\theta_n^{\,i} : \theta_m^{\,j} : 1) \mid 
\mbox{$i = 1, \dots , n$ and $j = 1, \dots , m$} \}\,.
\]
Then $\srees{I_H}$ is finitely generated.
\end{Example}

\section{Huneke's condition}

Let $K$ be a field and let $I$ be a proper homogeneous ideal
of the polynomial ring $S = K[x, y, z]$ which we regard as an
$\nn$-graded ring setting the degrees of $x$, $y$ and $z$ to
suitable positive integers.
We assume that $S / I$ is a $1$-dimensional reduced ring.
Let $\gm = (x, y, z)$ and $R = S_\gm$.
The following result can be proved by the same argument
developed in the proofs of \cite[Proposition 2.1 and Lemma 2.2]{KN}
replacing $x$ with $u$.
\begin{Theorem}\label{2.1}
Suppose $0 < r_i \in \nn$ and $\xi_i \in I^{(r_i)}$ for $i = 1, 2$.
Let us take a homogeneous element $u$ of $S$ so that
$uS + I$ is $\gm$-primary.
Then we have
\[
\length{R}{R / (u, \xi_1, \xi_2)R} \geq r_1r_2\cdot\length{S}{S / uS + I}
\]
and the following conditions are equivalent.
\begin{itemize}
\item[{\rm (1)}]
$\length{R}{R / (u, \xi_1, \xi_2)R} = r_1r_2\cdot\length{S}{S / uS + I}\,.$
\item[{\rm (2)}]
$\xi_1$ and $\xi_2$ satisfy Huneke's condition on $I$.
\end{itemize}
\end{Theorem}
As is described in Theorem \ref{1.3},
the finite generation of $\srees{I}$ can be characterized by
the existence of elements satisfying Huneke's condition on $I$.
Here, let us verify that Proposition \ref{1.2} follows from the equivalence
of the conditions (1) and (2) of Theorem \ref{2.1}.
In the rest of this paper, we assume $\deg x = \deg y = \deg z = 1$.
Suppose $\xi_i \in [I^{(r_i)}]_{d_i}$ for $i = 1, 2$,
where $0 < r_i, d_i \in \nn$.
If $u$ is a linear form in $S$ such that
$\length{R}{R / (u, \xi_1, \xi_2)R} < \infty$,
then $u, \xi_1, \xi_2$ is an $S$-regular sequence consisting of 
homogeneous polynomials of degrees $1, d_1, d_2$, respectively, and so
\[
\length{R}{R / (u, \xi_1, \xi_2)R} = \length{S}{S / (u, \xi_1, \xi_2)} = d_1d_2\,.
\]
On the other hand,
if $u$ is a linear form of $S$ whose
image in the local ring $R / IR$ generates
a reduction of the maximal ideal, we have
\[
\length{S}{S / uS + I} = \length{R}{R / uR + IR} =
 \mult{uR}{R / IR} = \mult{\gm}{R / IR} = \mult{}{S / I}\,.
\]
Consequently, if we choose a general linear form of $x$, $y$ and $z$
as $u$ of Theorem \ref{2.1},
the equality of (1) holds if and only if $d_1d_2 = r_1r_2\cdot\mult{}{S / I}$.
Thus we get Proposition \ref{1.2}.

In order to explain how to use Proposition \ref{1.2} and Theorem \ref{1.3},
let us verify the following well known example.

\begin{Example}\label{2.2}
Let $H$ be a set of of distinct $3$ points
$P_1, P_2, P_3 \in \pp^2$.
Then $\srees{I_H}$ is finitely generated.
\end{Example}

\noindent
{\it Proof}.\,\,
For $i \in \{ 1, 2, 3 \}$,
we take a linear form $f_i$ of $x$, $y$ and $z$
which defines the line going through $P_i$ and $P_{i + 1}$,
where $P_{i + 1}$ denotes $P_1$ for $i = 3$.
We set
\[
\xi_1 = f_1f_2f_3
\hspace{3ex}\mbox{and}\hspace{3ex}
\xi_2 = f_1f_2 + f_2f_3 + f_3f_1\,.
\]
Because $I_{P_1} = (f_1, f_2)$, $I_{P_2} = (f_2, f_3)$ and 
$I_{P_3} = (f_3, f_1)$,
it follows that
\[
\Min{}{(\xi_1, \xi_2)} = \{ I_{P_1}\,, I_{P_2}\,, I_{P_3} \}\,,
\]
and so $\height{}{(\xi_1, \xi_2)} = 2$.
On the other hand, as $f_i \in I_{P_i} \cap I_{P_{i + 1}}$
for any $i \in \{ 1, 2, 3 \}$,
we see
\[
\xi_1 \in I_{P_1}^{\,2} \cap I_{P_2}^{\,2} \cap I_{P_3}^{\,2} = I_H^{\,(2)}\,,
\]
and so $\xi_1 \in [ I_H^{\,(2)} ]_3$.
Similarly, we get $\xi_2 \in [ I_H ]_2$\,.
Because
\[
\frac{3}{2}\cdot\frac{2}{1} = 3 = \sharp H = \mult{}{S / I_H}\,,
\]
$\xi_1$ and $\xi_2$ satisfy Huneke's condition on $I_H$ by Proposition \ref{1.2}.
Therefore $\srees{I_H}$ is finitely generated by Theorem \ref{1.3}.
\qed

\section{An alternative proof of Theorem \ref{1.4}}

In the rest of this paper,
$K$ is an algebraically closed field and
the grading of $S = K[x, y, z]$ is ordinary.
We put $\gm = (x, y, z)$.
As is well known,
\[
\{\, \gp \in \mathrm{Spec}\, S \mid
\mbox{$\gp$ is homogeneous and $\dim S / \gp = 1$} \,\} =
\{\, I_P \mid P \in \pp^2 \,\}\,.
\]
For any $P \in \pp^2$,
we denote the localization of $S$ at $I_P$ and its maximal ideal
by $S_P$ and $\gm_P$, respectively.
Let $f$ and $g$ be non-zero homogeneous polynomials of $S$
such that $\deg f = m > 0$ and $\deg g = n > 0$.
We set
\[
H_{f, g} = \{\, P \in \pp^2 \mid (f, g) \subseteq I_P \,\}\,.
\]
Let us begin by verifying the following two lemmas,
which may be well known.

\begin{Lemma}\label{3.1}
The following conditions are equivalent.
\begin{itemize}
\item[{\rm (1)}]
$\dim S / (f, g) = 1$\,.
\item[{\rm (2)}]
$\Min{}{(f, g)} = \{\, I_P \mid P \in H_{f, g} \,\}$\,.
\item[{\rm (3)}]
$H_{f, g}$ is a finite set.
\end{itemize}
When this is the case,
$S / (f, g)$ is a Cohen-Macaulay ring.
\end{Lemma}

\noindent
{\it Proof}.\,\,
(1) $\Rightarrow$ (2)\,\,
Suppose $\dim S / (f, g) = 1$.
Let us take any $\gp \in \Min{}{(f, g)}$.
Then $\gp \subsetneq \gm$,
and so $0 < \dim S / \gp \leq \dim S / (f, g) = 1$.
Consequently, $\gp$ is a homogeneous ideal with $\dim S / \gp = 1$,
which means that $\gp = I_P$ for some $P \in \pp^2$.
Conversely, if $P \in H_{f, g}$,
we obviously have $I_P \in \Min{}{(f, g)}$.

(2) $\Rightarrow$ (3)\,\,
This implication holds since $\Min{}{(f, g)}$ is a finite subset of $\mathrm{Spec}\,S$.

(3) $\Rightarrow$ (1)\,\,
Suppose that $H_{f, g}$ is finite.
If $\height{}{(f, g)} = 1$,
there exists $h \in S$ such that $(f, g) \subseteq hS$,
which is impossible since there exist infinitely many $P \in \pp^2$
such that $h \in I_P$.
Thus we see $\height{}{(f, g)} = 2$,
and so $\dim S / (f, g) = 1$.
Then, as $f, g$ is an $S$-regular sequence,
$S / (f, g)$ is a Cohen-Macaulay ring.
\qed

\begin{Lemma}\label{3.2}
The following conditions are equivalent.
\begin{itemize}
\item[{\rm (1)}]
$S / (f, g)$ is a $1$-dimensional reduced ring.
\item[{\rm (2)}]
$\sharp H_{f, g} = mn$\,.
\item[{\rm (3)}]
$\dim S / (f, g) = 1$ and $\sharp H_{f, g} \geq mn$\,.
\item[{\rm (4)}]
$I_{H_{f, g}} = (f, g)$\,.
\end{itemize}
When this is the case, we have
$\gm_P = (f, g)S_P$ for any $P \in H_{f, g}$ and
$I_{H_{f, g}}^{\,\,(r)} = (f, g)^r$ for any $r \in \zz$.
\end{Lemma}

\noindent
{\it Proof}.\,\,
(1) $\Rightarrow$ (2)\,\,
Suppose that $S / (f, g)$ is a $1$-dimensional reduced ring.
Because $\dim S / (f, g) = 1$,
we have $\Min{}{(f, g)} = \{\, I_P \mid P \in H_{f, g} \,\}$
by Lemma \ref{3.1}.
Then, for any $P \in H_{f, g}$,
it follows that $S_P / (f, g)S_P$ is a field since
$S / (f, g)$ satisfies Serr's condition $(\mathrm{R}_0)$,
which means $\gm_P = (f, g)S_P$.
Here, let us choose a linear form $u \in S$ generally
so that its image in the Cohen-Macaulay local ring $R / (f, g)R$ generates a reduction
of the maximal ideal.
Then $u, f, g$ is a maximal $R$-regular sequence consisting of
homogeneous polynomials of degrees $1, m, n$, respectively,
and we have
\[
\mult{\gm}{R / (f, g)R} = \mult{uR}{R / (f, g)R)}
= \length{R}{R / (u, f, g)R} = \length{S}{S / (u, f, g)} = mn\,.
\]
On the other hand, by the additive formula of multiplicity, we have
\[
\mult{\gm}{R / (f, g)R} = 
\sum_{P \in H_{f, g}} \length{S_P}{S_P / \gm_P}\mult{\gm_P}{S_P / I_P}
= \sharp H_{f, g}\,.
\]
Thus we see that the condition (2) is satisfied.

(2) $\Rightarrow$ (3)\,\,
We get this implication by (3) $\Rightarrow$ (1) of Lemma \ref{3.1}.

(3) $\Rightarrow$ (4)\,\,
Suppose $\dim S / (f, g) = 1$ and $\sharp H_{f, g} \geq mn$.
Again, let us take a linear form $u \in S$ generally,
then we have
\begin{equation*}
\begin{split}
\mult{}{S / I_{H_{f, g}}} &= \mult{\gm}{R / (I_{H_{f, g}})R} = \mult{uR}{R / (I_{H_{f, g}})R}  \\
&\hspace{10ex}
= \length{R}{R / uR + (I_{H_{f, g}})R} =
\length{S}{uS + I_{H_{f, g}}}\,.
\end{split}
\end{equation*}
On the other hand, we have
\[
\mult{}{S / I_{H_{f, g}}} = \sharp H_{f, g} \geq mn = \length{S}{S / (u, f, g)}\,.
\]
Consequently, we get
\[
\length{S}{S / uS + I_{H_{f, g}}} \geq \length{S}{S / (u, f, g)}\,.
\]
However, as the inclusion $I_{H_{f, g}} \supseteq (f, g)$ holds obviously,
it follows that the both sides of the above inequality are equal, 
and so
$uS + I_{H_{f, g}} = (u, f, g)$.
Then
\begin{eqnarray*}
I_{H_{f, g}} & = & (u, f, g) \cap I_{H_{f, g}} \\
 & = & (f, g) + uS \cap I_{H_{f, g}} \\
 & = & (f, g) + u\cdot I_{H_{f, g}}\,.
\end{eqnarray*}
Therefore, by Nakayama's lemma, 
we see $I_{H_{f, g}} = (f, g)$.

(4) $\Rightarrow$ (1)\,\,
This implication is obvious.

Finally, we show $I_{H_{f, g}}^{\,\,(r)} = (f, g)^r$
for any $r \in \zz$ when the equivalent conditions (1) - (4) are satisfied.
Of course, we may assume $r > 0$.
Because $I_{H_{f, g}}^{\,\,(r)} \supseteq (f, g)^r$ holds obviously,
it is enough to show $I_{H_{f, g}}^{\,\,(r)}S_\gp = (f, g)^rS_\gp$,
where $\gp$ is any associated prime ideal of $S / (f, g)^r$.
In fact, as $S / (f, g)^r$ is a $1$-dimensional Cohen-Macaulay ring,
we have $\gp \in \Min{}{(f, g)}$,
and so there exists $P \in H_{f, g}$ such that $\gp = I_P$.
Then, $\gm_P = (f, g)S_P$ as is proved in the proof of (1) $\Rightarrow$ (2).
Hence we have
\[
I_{H_{f, g}}^{\,\,(r)}S_P = I_P^{\,r}S_P = \gm_P^{\,r} = (f, g)^rS_P\,,
\]
and so the proof is complete as $S_\gp = S_P$.
\qed

Now, we are ready to give an alternative proof for Theorem \ref{1.4}
using Huneke's criterion.
In the rest of this section,
let $n$ be a positive integer which is not a multiple of
the characteristic of $K$.
We take a primitive $n$-th root $\theta$ of unity, and set
\[
H = \{ (1 : 0 : 0), (0 : 1 : 0), (0 : 0 : 1) \} \cup
\{ P_{ij} \mid i, j = 1, \dots , n \} \subset \pp^2\,,
\]
where $P_{ij} = (\theta^i : \theta^j : 1)$.
Let
\[
f = y^n - z^n \,,\,\, g = z^n - x^n
\hspace{1ex}\mbox{and}\hspace{1ex} 
h = x^n - y^n\,.
\]
Then, as $f + g + h = 0$, we have
\begin{equation}\label{eq3.1}
(f, g) = (g, h) = (h, f)
\hspace{1ex}\mbox{and}\hspace{1ex}
H_{f, g} = H_{g, h} = H_{h, f}\,.
\end{equation}
Moreover, it is easy to see that
\begin{equation}\label{eq3.2}
\mbox{
$f$, $g$ and $h$ are elements of $I_{P_{ij}}$
for any $i, j = 1, \dots , n$\,,
}
\end{equation}
which means $\{ P_{ij} \}_{i, j } \subseteq H_{f, g}$.
Because $\dim S / (f, g) = 1$ and
$\sharp\{ P_{ij} \}_{i, j } = n^2$,
by Lemma \ref{3.2} we see
\begin{equation}\label{eq3.3}
\mbox{
$H_{f, g} = \{ P_{ij} \}_{i, j}$\,,
$I_{H_{f, g}} = (f, g)$ and
$\gm_{P_{ij}} = (f, g)S_{P_{ij}}$ for any $i, j $\,.
}
\end{equation}
Because $I_{(1\,:\,0\,:\,0)} = (y, z)$, $I_{(0\,:\,1\,:\,0)} = (z, x)$ and
$I_{(0\,:\,0\,:\,1)} = (x, y)$,
we get the following assertions by 
Lemma \ref{3.2}, \eqref{eq3.1}, \eqref{3.2} and \eqref{eq3.3}.
\begin{equation}\label{eq3.4}
\mbox{
$I_H^{\,(r)} = (y, z)^r \cap (z, x)^r \cap (x, y)^r \cap (f, g)^r$
for any $r \in \zz$\,.
}
\end{equation}
\begin{equation}\label{eq3.5}
\mbox{
$xf$, $yg$ and $zh$ are elements of $I_H$\,.
}
\end{equation}

If $n = 1$ or $2$, then $\sharp H = 4$ or $7$,
and so $\srees{I_H}$ is finitely generated as is mentioned in Introduction.
Hence, we may assume $n \geq 3$.

First, let us consider the case where $n = 3$.
In this case, we set
\[
\xi_1 = fgh
\hspace{2ex}\mbox{and}\hspace{2ex}
\xi_2 = xf\cdot yg + yg \cdot zh + zh\cdot xf\,.
\]
By \eqref{eq3.4} and \eqref{eq3.5},
we have $\xi_1 \in [ I_H^{\,(3)} ]_9$ and 
$\xi_2 \in [ I_H^{\,2} ]_8 \subseteq [ I_H^{\,(2)} ]_8$.
Let $\gp$ be any prime ideal of $S$ containing $\xi_1$ and $\xi_2$.
Because $\xi_1 \in \gp$,
one of $f$, $g$ and $h$ belongs to $\gp$.
If $f \in \gp$,
then $yg\cdot zh \in\gp$ as $\xi_2 \in \gp$,
and so $\height{}{\gp} \geq 2$
as $\gp$ includes one of $(f, y)$, $(f, z)$ or $(f, g)$
($= (f, h)$).
Similarly, we get $\height{}{\gp} \geq 2$ if $g \in \gp$ or $h \in \gp$.
Consequently, we have $\height{}{(\xi_1, \xi_2)} = 2$.
Hence, by Proposition \ref{1.2}
it follows that $\xi_1$ and $\xi_2$ satisfy {\bf HC} on $I_H$ since
\[
\frac{9}{3}\cdot\frac{8}{2} = 12 = \sharp H = \mult{}{S / I_H}\,.
\]
Therefore $\srees{I_H}$ is finitely generated by Theorem \ref{1.3}.

In the rest of this section,
we assume $n \geq 4$.
In this case, taking an element $\alpha \in K$ so that $\alpha \neq 0, 1$,
we set
\begin{equation*}
\begin{split}
&\xi_1 = fgh\cdot(\alpha f + g)^{n - 3} \hspace{2ex}\mbox{and} \\
&\xi_2 = (xf)^2\cdot(yg)^{n - 2} + (yg)^2\cdot(zh)^{n - 2} 
  + (zh)^2\cdot(xf)^{n - 2} + f^{n - 2}gh\,.
\end{split}
\end{equation*}
Let us notice that $\xi_2$ is not homogeneous although so is $\xi_1$.
By \eqref{eq3.4} and \eqref{eq3.5} we can easily verify that
\begin{equation}\label{3.6}
\mbox{
$\xi_1$ and $\xi_2$ belongs to $I_H^{\,(n)}$\,.
}
\end{equation}
We aim to
show that $\xi_1$ and $\xi_2$ satisfy Huneke's condition on $I_H$.

First, let us verify $I_HR = \sqrt{(\xi_1, \xi_2)R}$, where $R = S_\gm$.
As is noticed in Definition \ref{1.1},
the crucial point is to prove that
the right side includes the left side.
For that purpose,
it is enough to see that the following assertion is true by \eqref{eq3.4}.
\begin{Claim}\label{3.3}
Let $\gp$ be a prime ideal of $S$ such that
$(\xi_1, \xi_2) \subseteq \gp \subseteq \gm$.
Then $\gp$ includes one of 
$(x, y)$, $(y, z)$, $(z, x)$ or $(f, g)$.
\end{Claim}
In fact, as $\xi_1 \in \gp$,
one of $f$, $g$, $h$ or $\alpha f + g$ belongs to $\gp$.
If $f \in \gp$,
then $(yg)^2\cdot(zh)^{n - 2} \in \gp$ as $\xi_2 \in \gp$,
and so $\gp$ includes $(y, z)$ or $(f, g)$ since
$\sqrt{(f, y)} = \sqrt{(f, z)} = (y, z)$ and $(f, g) = (h, f)$ by \eqref{eq3.1}.
Similarly, we see that $\gp$ includes one of $(x, y)$, $(z, x)$ or $(f, g)$
if $g \in \gp$ or $h \in \gp$.
So, let us consider the case where $\alpha f + g \in \gp$.
Then, as 
\[
g \equiv -\alpha f \hspace{1ex}\mbox{mod}\hspace{1ex} \gp
\hspace{3ex}\mbox{and}\hspace{3ex}
h = -(f + g) \equiv (\alpha - 1)f \hspace{1ex}\mbox{mod}\hspace{1ex} \gp\,,
\]
it follows that
\begin{equation*}
\begin{split}
\xi_2 &\equiv f^n\eta \hspace{1ex}\mbox{mod}\hspace{1ex}\gp\,,\,\mbox{where}  \\
&\hspace{1ex}\eta =
(-\alpha)^{n - 2}x^2y^{n - 2} + \alpha^2(\alpha - 1)^{n - 2}y^2z^{n - 2} +
(\alpha - 1)^2z^2x^{n - 2} - \alpha(\alpha - 1)\,.
\end{split}
\end{equation*}
Because $\alpha(\alpha - 1) \neq 0$, we have $\eta \not\in \gm$,
which means $\eta \not\in \gp$.
Hence, we get $f \in \gp$ as $\xi_2 \in \gp$,
and so $\gp$ includes $(\alpha f + g, f) = (f, g)$.
Thus we have seen Claim \ref{3.3}.

Next, we verify $\gr{S_P}_+ = \sqrt{(\xi_1t^n, \xi_2t^n)\gr{S_P}}$
for any $P \in H$.
Again, the crucial point is to prove that
the right side includes the left side,
which is deduced from the next assertion.
\begin{Claim}\label{3.4}
Let $P \in H$ and $\calP$ be a prime ideal of $\gr{S_P}$
containing $(\xi_1t^n, \xi_2t^n)\gr{S_P}$.
Then we have $\calP = \gr{S_P}_+$.
\end{Claim}
If $P \in H$ and $\eta \in \gm_P^{\,r} \cap S$ for $r \in \nn$,
we denote by $\overline{\eta t^r}$ the image of 
$(\eta / 1)t^r \in \gm_P^{\,r}t^r$ under the homomorphism
$\rees{S_P} \rightarrow \gr{S_P}$.

Let us start the proof of Claim \ref{3.4} 
with checking the case where $P \in H_{f, g}$.
In this case, we have
$\gm_P = (f, g)S_P = (g, h)S_P = (h, f)S_P$
by \eqref{3.1} and Lemma \ref{3.2},
and none of $x$, $y$ and $z$ belongs to $I_P$,
which means that $\overline{x}$,
$\overline{y}$ and $\overline{z}$ are unites of $\gr{S_P}$.
We set
\[
U = \overline{ft}\,,\,\, V = \overline{gt}\quad \mbox{and}\quad W = \overline{ht}\,.
\]
Then we have $U + V + W = 0$ and
$\gr{S_P}_+$ is generated by any two elements of $\{ U, V, W \}$
as an ideal of $\gr{S_P}$.
Moreover, we have the equalities
\begin{eqnarray*}
\overline{\xi_1t^n} & = & UVW(\overline{\alpha}\cdot U + V)^{n -3} \quad \mbox{and} \\
\overline{\xi_2t^n} & = & \overline{xy^{n - 2}}\cdot U^2V^{n - 2}
  + \overline{y^2z^{n - 2}}\cdot V^2W^{n - 2}
  + \overline{z^2x^{n - 2}}\cdot W^2U^{n - 2} + U^{n - 2}VW
\end{eqnarray*}
in $\gr{S_P}$.
Because $\overline{\xi_1t^n} \in \calP$,
one of $U$, $V$, $W$ and $\overline{\alpha}\cdot U + V$ belongs to $\calP$.
If $U \in \calP$,
then $\overline{y^2z^{n - 2}}\cdot V^2W^{n - 2} \in \calP$
as $\overline{\xi_2t^n} \in \calP$,
and so $\calP$ includes $\gr{S_P}_+ = (U, V) = (U, W)$ as
$\overline{y^2z^{n - 2}}$ is a unit of $\gr{S_P}$.
Similarly, we can see that $\calP$ includes $\gr{S_P}_+$ if
$V \in \calP$ or $W \in \calP$.
So, let us consider the case where $\overline{\alpha}\cdot U + V \in \calP$.
Then, as
\[
V \equiv -\overline{\alpha}\cdot U\,\mathrm{mod}\,\calP
\quad\mbox{and}\quad
W = -(U + V) \equiv \overline{\alpha - 1}\cdot U \,\mathrm{mod} \, \calP\,,
\]
it follows that
\[
\overline{\xi_2t^n} \equiv \overline{\eta}\cdot U^n \, \mathrm{mod} \, \calP\,,
\]
where $\eta$ is the element stated in the proof of Claim \ref{3.3}.
Because $\eta \not\in \gm$,
$\overline{\eta}$ is a unit of $\gr{S_P}$.
Hence we get $U \in \calP$ as $\overline{\xi_2t^n} \in \calP$.
Thus we see that $\calP$ includes $\overline{\alpha}\cdot U + V$ and $U$,
which means $\calP = \gr{S_P}_+$.

Next, let us consider the case where $P = (0 : 0 : 1)$.
Then, none of $f$, $g$ and $z$ belongs to $I_P = (x, y)$.
Moreover, $\alpha f + g \not\in I_P$ since
$\alpha f + g \equiv (1 - \alpha)z^n$ $\mathrm{mod}$ $I_P$ and $\alpha \neq 1$.
Hence $\overline{f}$, $\overline{g}$, $\overline{z}$ and 
$\overline{\alpha f + g}$ are units of $\gr{S_P}$.
On the other hand, $\gm_P = (x, y)S_P$.
So, we set
\[
X = \overline{xt} \quad\mbox{and}\quad Y = \overline{yt}\,.
\]
Then $\gr{S_P} = (X, Y)$, and we have the equalities
\begin{eqnarray*}
\overline{\xi_1t^n} & = & 
    \overline{fg(\alpha f + g)^{n - 3}}\cdot(X^n - Y^n) \quad \mbox{and} \\
\overline{\xi_2t^n} & = & 
    \overline{f^2g^{n - 2}}\cdot X^2Y^{n - 2} + \overline{f^{n - 2}g}\cdot(X^n - Y^n)
\end{eqnarray*}
in $\gr{S_P}$,
where the second equality holds since $(yg)^2\cdot(zh)^{n - 2}$ and
$(zh)^2\cdot(xf)^{n - 2}$ are included in $I_P^{\,n + 1}$.
Because $\overline{\xi_1t^n} \in \calP$ and
$\overline{fg(\alpha f + g)^{n - 2}}$ is a unit of $\gr{S_P}$,
we have $X^n - Y^n \in \calP$.
Then we get $X^2Y^{n - 2} \in \calP$ since 
$\overline{\xi_2t^n} \in \calP$ and $\overline{f^2g^{n - 2}}$ is a unit of $\gr{S_P}$.
Consequently, we see that $X$ and $Y$ belong to $\calP$,
which means $\calP = \gr{S_P}_+$.

If $P = (0 : 1 : 0)$, we can prove $\calP = \gr{S_P}_+$ similarly as the above case.

Finally, we suppose $P = (1 : 0 : 0)$.
In this case, $I_P = (y, z)$ and 
$\overline{g}$, $\overline{h}$, $\overline{x}$ and $\overline{\alpha f + g}$
are units in $\gr{S_P}$.
On the other hand, $\gm_P = (y, z)S_P$.
So, we set
\[
Y = \overline{yt} \quad\mbox{and}\quad Z = \overline{zt}\,.
\]
Then $\gr{S_P}_+ = (Y, Z)$, and we have the equalities
\begin{eqnarray*}
\overline{\xi_1t^n} & = &
  \overline{gh(\alpha f + g)^{n - 3}}\cdot(Y^n - Z^n) \quad\mbox{and} \\
\overline{\xi_2t^n} & = & \overline{g^2h^{n - 2}}\cdot Y^2Z^{n - 2}
\end{eqnarray*}
in $\gr{S_P}$, where the second equality holds since
$(xf)\cdot(yg)^{n - 2}$, $(zh)^2\cdot(xf)^{n - 2}$ and $f^{n - 2}gh$
are included in $I_P^{\,n + 1}$.
Because $\overline{\xi_1t^n} \in \calP$
and $\overline{gh(\alpha f + g)^{n - 3}}$ is a unit of $\gr{S_P}$,
we get $Y^n - Z^n \in \calP$.
On the other hand, we get $Y^2Z^{n - 2} \in \calP$ since
$\overline{\xi_2t^n} \in \calP$ and $\overline{g^2h^{n - 2}}$ is a unit of $\gr{S_P}$.
Consequently, we see that $Y$ and $Z$ belong to $\calP$,
which means $\calP = \gr{S_P}_+$.
Thus the proof of Theorem \ref{1.4} is complete.
\begin{Remark}\label{3.5}
Suppose $n \geq 4$.
Then, $I_H^{\,(n)}$ has no reduction generated by
two homogeneous polynomials by \cite[Proposition 5.1]{NS}.
However, by the argument stated in the proof of \cite[Theorem 2.5]{KN},
we can prove that $(\xi_1, \xi_2)R$ is a reduction of $I_H^{\,(n)}R$.
\end{Remark}

\section{Proof of Theorem 1.5}

In this section, let $f$ and $g$ be homogeneous polynomials of $S$
having positive degrees $m$ and $n$, respectively.
We assume $f \in I_A^{\,m}$ and $g \in I_B^{\,n}$,
where $A$ and $B$ are points of $\pp^2$.
Let us take linear forms $u, v \in S$ so that
$I_A = (u, v)$.
Because $f \in I_A^{\,m}$, we can express
\[
f = \sum_{j = 0}^m a_ju^jv^{m - j} \quad (\, a_j \in S \,)\,.
\]
However, as $f$ is a homogeneous polynomial of degree $m$,
we can choose $a_0, a_1, \dots, a_m$ from $K$.
Then
\[
\frac{f}{v^m} = \sum_{j = 0}^m a_j \cdot(\,\frac{u}{v}\,)^{\,j} \in K[\, \frac{u}{v} \,]\,.
\]
Because $K$ is algebraically closed, we can express
\[
\frac{f}{v^m} = \prod_{i = 1}^m (\alpha_i\cdot\frac{u}{v} - \beta_i)
\quad (\,\alpha_i, \beta_i \in K \,)\,.
\]
Then, setting $f_i = \alpha_i u - \beta_i v \in [ I_A ]_1$
for $i = 1, 2, \dots, m$, we have
\[
f = f_1f_2\cdots f_m\,.
\]
Similarly, there exist linear forms $g_1, g_2, \dots, g_n \in [ I_B ]_1$ such that
\[
g = g_1g_2\cdots g_n\,.
\]

In the rest of this section, we assume
\[
A \neq B\,,\hspace{1ex} f \not\in I_B 
\hspace{1ex}\mbox{and}\hspace{1ex} g \not\in I_A\,.
\]
Then we have
\begin{equation}\label{eq4.1}
A, B \not\in H_{f, g}\,.
\end{equation}
Moreover, for any $i = 1, \dots, m$ and $j = 1, \dots, n$,
we have $f_i \not\in I_B$ and $g_j \not\in I_A$,
and so $f_i \not\sim g_j$,
which means that $f_i$ and $g_j$ define distinct two lines in $\pp^2$
intersecting at the point $P_{ij}$ with $I_{P_{ij}} = (f_i, g_j)$.
Of course $P_{ij} \in H_{f, g}$ for any $i. j$.

Let us assume furthermore that
$S / (f, g)$ is a $1$-dimensional reduced ring.
Then the following assertions hold by Lemma \ref{3.2}.
\begin{equation}\label{eq4.2}
\sharp H_{f, g} = mn\,.
\end{equation}
\begin{equation}\label{eq4.3}
\mbox{$\gm_P = (f, g)S_P$ for any $P \in H_{f, g}$\,.}
\end{equation}
\begin{equation}\label{eq4.4}
\mbox{$I_{H_{f, g}}^{\,(r)} = (f, g)^r$ for any $r \in \zz$\,.}
\end{equation}
Moreover, we have $f_i \not\sim f_k$ if $i \neq k$ and
$g_j \not\sim g_\ell$ if $j \neq \ell$.
Here we suppose $P_{ij} = P_{k\ell}$.
Then $f_k \in I_{P_{k\ell}} = I_{P_{ij}}$.
Hence, if $i \neq k$, we have $I_{P_{ij}} = (f_i, f_k) = I_A$
as $f_i \not\sim f_k$,
which contradicts to \eqref{eq4.1}.
Thus we get $i = k$.
Similarly, we get also $j = \ell$.
Consequently, we see $P_{ij} \neq P_{k\ell}$
if $i \neq k$ or $j \neq \ell$,
and so $\sharp\{ P_{ij} \}_{i, j} = mn$.
Hence the following assertion is deduced by \eqref{eq4.2}.
\begin{equation}\label{eq4.5}
H_{f, g} = \{ P_{ij} \mid 
\mbox{$i = 1, \dots, m$ and $j = i, \dots, n$} \}\,.
\end{equation}

Let $h$ be a linear form in $S$ defining a line going through $A$ and $B$,
i.e., $h \in [ I_A \cap I_B ]_1$.
For any $i = 1, \dots, m$,
we have $f_i \not\sim h$ since $f_i \not\in I_B$ and $h \in I_B$.
Hence we see
\begin{equation}\label{eq4.6}
\mbox{
$I_A = (f_i, h)$ for any $i = 1, \dots, m$\,.
}
\end{equation}
As a consequence, we get
\begin{equation}\label{eq4.7}
(f, h) \subseteq \gp \in \mathrm{Spec}\,S
\hspace{1ex}\Rightarrow\hspace{1ex}
I_A \subseteq \gp\,.
\end{equation}
The following two assertions can be verified similarly as \eqref{eq4.6} and \eqref{eq4.7}.
\begin{equation}\label{eq4.8}
\mbox{
$I_B = (g_j, h)$ for any $j = 1, \dots, n$\,.
}
\end{equation}
\begin{equation}\label{eq4.9}
(g, h) \subseteq \gp \in \mathrm{Spec}\,S
\hspace{1ex}\Rightarrow\hspace{1ex}
I_B \subseteq \gp\,.
\end{equation}
Let us take any $P \in H_{f, g}$.
If $h \in I_P$,
then $I_A = I_P$ by \eqref{eq4.7},
which contradicts to \eqref{eq4.1}.
Hence we have
\begin{equation}\label{eq4.10}
\mbox{
$h \not\in I_P$ for any $P \in H_{f, g}$\,.
}
\end{equation}

We set $H = \{ A, B \} \cup H_{f, g}$.
By \eqref{eq4.1} and \eqref{eq4.4}, we have
\begin{equation}\label{eq4.11}
\mbox{
$I_H^{\,(r)} = I_A^{\,r} \cap I_B^{\,r} \cap (f, g)^r$ for any $r \in \zz$\,.
}
\end{equation}
Similarly as in Section 3,
if $P \in H$ and $\eta \in \gm_P^{\,r} \cap S$ for $r \in \nn$,
we denote by $\overline{\eta t^r}$ the image of
$(\eta / 1)t^r \in \gm_P^{\,r}t^r$ under the homomorphism
$\rees{S_P} \rightarrow \gr{S_P}$.
Here we want to show the following assertion.
\begin{equation}\label{eq4.12}
\mbox{
$\overline{ft^m}$, $\overline{ht}$ is an sop for $\gr{S_A}$\,.
}
\end{equation}
It is enough to show that
$\gr{S_A}_+$ is the unique prime ideal of 
$\gr{S_A}$ containing $\overline{ft^m}$ and $\overline{ht}$.
So, let us take any $\calP \in \mathrm{Spec}\,\gr{S_A}$
containing $\overline{ft^m}$ and $\overline{ht}$.
Because the factorization
\[
\overline{ft^m} = \prod_{i = 1}^m \overline{f_i t}
\]
holds in $\gr{S_A}$,
we can choose $i = 1, \dots, m$ so that
$\overline{f_i t} \in \calP$.
Then, we have $\calP = \gr{S_A}_+$ since $\gm_A = (f_i, h)S_A$ by \eqref{eq4.6}.
Similarly, the following assertion holds.
\begin{equation}\label{eq4.13}
\mbox{
$\overline{gt^n}$, $\overline{ht}$ is an sop for $\gr{S_B}$\,.
}
\end{equation}

Now, we are ready to prove Theorem \ref{1.5}.
If $m = 1$ or $n = 1$,
then all the points of $H$ except for just one point lie on a line,
and so the Cox ring $\Delta_H$ is finitely generated
by the result due to Testa, Varilly-Alvarado and Verasco
(The case (ii) stated in Introduction can be applied).
So, in the rest, we assume $m \geq 2$ and $n \geq 2$.
We set
\[
\xi_1 = f^ng^m(f + g)^{mn - m - n}
\hspace{2ex}\mbox{and}\hspace{2ex}
\xi_2 = fg + (f + g)^2h^2\,.
\]
Because $f \in I_A^{\,m}$, $g \in I_B^{\,n}$ and $h \in I_A \cap I_B$,
we have
\[
\xi_1 \in I_H^{\,(mn)}
\quad\mbox{and}\quad
\xi_2 \in I_H^{\,(2)}
\]
by \eqref{eq4.11}.
We aim to show that $\xi_1$ and $\xi_2$ satisfy Huneke's condition on $I_H$.

First, let us verify $I_H = \sqrt{(\xi_1, \xi_2)}$,
which implies $I_HR = \sqrt{(\xi_1, \xi_2)R}$.
For that purpose, 
it is enough to see that the following assertion is true by 
\eqref{eq4.7}, \eqref{eq4.9} and \eqref{eq4.11}.
\begin{Claim}\label{4.1}
Let $\gp$ be a prime ideal of $S$ containing
$\xi_1$ and $\xi_2$.
Then $\gp$ includes one of $(f, h)$, $(g, h)$ or $(f, g)$\,.
\end{Claim}
\noindent
In fact, as $\xi_1 \in \gp$,
one of $f$, $g$ or $f + g$ belongs to $\gp$.
If $f \in \gp$,
then $(f + g)^2h^2 \in \gp$ as $\xi_2 \in \gp$,
and so $\gp$ includes $(f, f + g) = (f, g)$ or $(f, h)$.
Similarly, we see that $\gp$ includes $(f, g)$ or $(g, h)$ if $g \in \gp$.
If $f + g \in \gp$,
then $fg \in \gp$ as $\xi_2 \in \gp$,
and so $\gp$ includes $(f, g)$ as $f \in \gp$ or $g \in \gp$.
Thus we have seen Claim \ref{4.1}.

Next, we verify $\gr{S_P}_+ = \sqrt{(\xi_1t^{mn}, \xi_2t^2)\gr{S_P}}$
for any $P \in H$,
which is deduced from the next assertion.
\begin{Claim}\label{4.2}
Let $P \in H$ and $\calP$ be a prime ideal of $\gr{S_P}$
containing $\overline{\xi_1t^{mn}}$ and $\overline{\xi_2t^2}$.
Then we have $\calP = \gr{S_P}_+$.
\end{Claim}
\noindent
Let us start the proof of the above assertion
with checking the case where $P \in H_{f, g}$.
In this case, we have $\gm_P = (f, g)S_P$ by \eqref{eq4.3} and
$\overline{h}$ is a unit of $\gr{S_P}$ by \eqref{eq4.10}.
We set
\[
U = \overline{ft}
\quad\mbox{and}\quad
V = \overline{gt}\,.
\]
Then $\gr{S_P}_+ = (U, V) = (U, U + V) = (U + V, V)$
and we have the equalities
\[
\overline{\xi_1t^{mn}} = U^nV^m(U + V)^{mn - m - n}
\quad\mbox{and}\quad
\overline{\xi_2t^2} = UV + (U + V)^2\cdot \overline{h^2}
\]
in $\gr{S_P}$.
Because $\overline{\xi_1t^{mn}} \in \calP$,
one of $U$, $V$ or $U + V$ belongs to $\calP$.
If $U \in \calP$,
then $U + V \in \calP$ since $\overline{\xi_2t^2} \in \calP$
and $\overline{h^2} = (\overline{h})^2$ is a unit,
and so $\calP = \gr{S_P}_+$.
Similarly, we see $\calP = \gr{S_P}_+$ if $V \in \calP$.
If $U + V \in \calP$,
then $UV \in \calP$ as $\overline{\xi_2t^2} \in \calP$,
and so $\calP = \gr{S_P}_+$ as $P$ contains $U$ or $V$.

Next, we consider the case where $P = A$.
Let us notice that the equalities
\[
\overline{\xi_1t^{mn}} = (\overline{ft^m})^n\cdot  \overline{g^m(f + g)^{mn - m - n}}
\quad\mbox{and}\quad
\overline{\xi_2t^2} = \overline{ft^2}\cdot\overline{g} + 
   \overline{(f + g)^2}\cdot (\overline{ht})^2
\]
hold in $\gr{S_A}$.
Because $I_A$ does not include $g$ and $f + g$,
it follows that
$\overline{g^m(f + g)^{mn - m - n}}$, $\overline{g}$ and
$\overline{(f + g)^2}$ are units in $\gr{S_A}$.
Hence we have $\overline{ft^m} \in \calP$ as $\overline{\xi_1t^{mn}} \in \calP$.
Then $\overline{ft^2}$ also belongs to $\calP$ 
as it vanishes if $m \geq 3$,
and so $\overline{ht} \in \calP$ as $\overline{\xi_2t^2} \in \calP$.
Therefore we see $\calP = \gr{S_A}_+$ by \eqref{eq4.12}.

Finally, the case where $P = B$ can be verified as above using \eqref{eq4.13}.
Thus we have seen Claim \ref{4.2}, and the proof of Theorem \ref{1.5} is complete.

\end{document}